# A Note on Hadamard Inequalities for the Product of the Convex Functions

Şahin Emrah Amrahov

*Computer Engineering Department of Engineering Faculty, Ankara University, 06100 Tandoğan, Ankara, Turkey*

*e-mail: emrah@eng.ankara.edu.tr*

*Abstract.* The main aim of the present note is to prove new Hadamard like integral inequalities for the product of the convex functions.

## 1. Introduction

Let $f$ be a real valued convex function defined on the interval $[a,b]$. Then

$$f\left(\frac{a+b}{2}\right) \leq \frac{1}{b-a}\int_a^b f(x)dx \leq \frac{f(a)+f(b)}{2} \tag{1}$$

is known as Hadamard's inequality for the convex function. [1]

If $f(x) = u(x)v(x)$ and $u(x), v(x)$ are convex functions then we have

$$u\left(\frac{a+b}{2}\right)v\left(\frac{a+b}{2}\right) \leq \frac{1}{b-a}\int_a^b u(x)v(x)dx \leq \frac{u(a)v(a)+u(b)v(b)}{2} \tag{2}$$

by Hadamard inequality.

Therefore by Cauchy-Scwartz inequality

$$\frac{u(a)v(a)+u(b)v(b)}{2} \leq \frac{\sqrt{u^2(a)+u^2(b)}\sqrt{v^2(a)+v^2(b)}}{2} \tag{3}$$

the inequality

$$\frac{1}{b-a}\int_a^b u(x)v(x)dx \leq \frac{\sqrt{u^2(a)+u^2(b)}\sqrt{v^2(a)+v^2(b)}}{2} \tag{4}$$

holds.

Note that if $u(x), v(x)$ are convex functions, then the function $f(x) = u(x)v(x)$ may not be convex function.

**Example1**. $u(x) = x^2, v(x) = (1-x)^2$ are convex functions defined on $[0,1]$, but the function $f(x) = u(x)v(x) = x^2(1-x)^2 = (x^2-x)^2$ is not convex function on this interval. In fact, the first and second derivatives $f'(x), f''(x)$ of $f(x)$ can be calculated by the formula

$$f'(x) = 2(x^2 - x)(2x - 1) \quad f''(x) = 2(2x-1)^2 + 4(x^2 - x)$$

Therefore we have $f''\left(\dfrac{1}{2}\right) = -1$.

It means that the function $f(x) = u(x)v(x)$ is not convex function on the interval $[0,1]$.

**Example2**. For the convex functions $u(x) = x^2, v(x) = (1-x)^2$ defined on $[0,1]$ we have

$$\frac{1}{b-a}\int_a^b u(x)v(x)dx = \int_0^1 (x^2 - x)^2 dx = \frac{1}{30}$$

$u(0)v(0) = 0$,
$u(1)v(1) = 0$.
Hence the Hadamard inequality

$$\frac{1}{b-a}\int_a^b u(x)v(x)dx = \int_0^1 (x^2 - x)^2 dx = \frac{1}{30} \leq \frac{u(0)v(0) + u(1)v(1)}{2} = 0$$

for the non-convex function $f(x) = u(x)v(x)$ does not hold.
On the other hand, since
$u^2(0) = 0$,
$u^2(1) = 1$,
$v^2(0) = 1$,
$v^2(1) = 0$
the inequality

$$\frac{1}{b-a}\int_a^b u(x)v(x)dx = \int_0^1 (x^2 - x)^2 dx = \frac{1}{30} \leq \frac{\sqrt{u^2(0) + u^2(1)}\sqrt{v^2(0) + v^2(1)}}{2} = \frac{1}{2}$$

holds. It means that although the function $f(x) = u(x)v(x)$ is no convex function we have that the inequality (4) is true for these functions.
Our aim is to investigate the inequality (4) when $f(x) = u(x)v(x)$ is non-convex function.

2. **Main results**
Our main result is the following theorem.
**Theorem.** Let $u(x)$ and $v(x)$ are nonnegative convex functions defined on the interval $[a,b]$. Then the inequality (4) holds.
To prove of the theorem we need the following lemma.
**Lemma**. Let $u(x)$ is a nonnegative convex function defined on the interval $[a,b]$. Then the function $u^2(x)$ is also convex function on the interval $[a,b]$.
**Proof**. For arbitrary $x, y \in [a,b]$ and $k \in [0,1]$ we have

$$(u(x)-u(y))^2 \geq 0$$

$$u^2(x) - 2u(x)u(y) + u^2(y) \geq 0$$

$$2u(x)u(y) \leq u^2(x) + u^2(y) \qquad (5)$$

Multiplying both sides of the inequality (5) by $k(1-k)$ we get

$$2k(1-k)u(x)u(y) \leq k(1-k)u^2(x) + k(1-k)u^2(y)$$

Therefore

$$2k(1-k)u(x)u(y) \leq (k-k^2)u^2(x) + [(1-k)-(1-k)^2]u^2(y)$$

$$2k(1-k)u(x)u(y) \leq ku^2(x) - k^2u^2(x) + (1-k)u^2(y) - (1-k)^2 u^2(y)$$

Hence

$$k^2 u^2(x) + 2k(1-k)u(x)u(y) + (1-k)^2 u^2(y) \leq ku^2(x) + (1-k)u^2(y)$$

$$[ku(x) + (1-k)u(y)]^2 \leq ku^2(x) + (1-k)u^2(y) \qquad (6)$$

Since $u(x)$ is a nonnegative convex function we have

$$u(kx + (1-k)y) \leq ku(x) + (1-k)u(y)$$

$$u^2(kx + (1-k)y) \leq [ku(x) + (1-k)u(y)]^2 \qquad (7)$$

From (6) and (7) we get

$$u^2(kx + (1-k)y) \leq ku^2(x) + (1-k)u^2(y) \qquad (8)$$

The inequality (8) proves that the function $u^2(x)$ is a convex function.

**Proof of the theorem.**

By the lemma the functions $u^2(x)$ and $v^2(x)$ are convex functions. By the Hadamard inequality for these functions we have

$$\frac{1}{b-a}\int_a^b u^2(x)dx \leq \frac{u^2(a) + u^2(b)}{2} \qquad (9)$$

$$\frac{1}{b-a}\int_a^b v^2(x)dx \leq \frac{v^2(a) + v^2(b)}{2} \qquad (10)$$

Multiplying the inequalities (9) and (10) we get

$$\frac{1}{(b-a)^2}\int_a^b u^2(x)dx \int_a^b v^2(x)dx \leq \frac{u^2(a) + u^2(b)}{2} \frac{v^2(a) + v^2(b)}{2} \qquad (11)$$

By Cauchy-Scwartz inequality

$$\left(\int_a^b u(x)v(x)dx\right)^2 \leq \int_a^b u^2(x)dx \int_a^b v^2(x)dx \qquad (12)$$

Hence by (11) and (12) we get

$$\frac{1}{(b-a)^2}\left(\int_a^b u(x)v(x)dx\right)^2 \leq \frac{u^2(a) + u^2(b)}{2} \frac{v^2(a) + v^2(b)}{2}$$

The last inequality means that the inequality (4) holds.

**Example 3.** Prove the inequality $\int_\pi^{2\pi} \frac{\sin x + 8}{x} dx \leq 2\sqrt{10}$.

Solution. Since the functions $u(x) = \sin x + 8$ and $v(x) = \dfrac{1}{x}$ are nonnegative convex functions on $[\pi, 2\pi]$ we have by inequality (4)

$$\frac{1}{\pi}\int_{\pi}^{2\pi}\frac{\sin x + 8}{x}dx \leq \frac{\sqrt{(\sin \pi + 8)^2 + (\sin 2\pi + 8)^2}\sqrt{\dfrac{1}{\pi^2} + \dfrac{1}{4\pi^2}}}{2}$$

Hence

$$\frac{1}{\pi}\int_{\pi}^{2\pi}\frac{\sin x + 8}{x}dx \leq \frac{\sqrt{8^2 + 8^2}\sqrt{5}}{4\pi}$$

Therefore

$$\int_{\pi}^{2\pi}\frac{\sin x + 8}{x}dx \leq 2\sqrt{10}$$

Note that $f(x) = \dfrac{\sin x + 8}{x}$ is not a convex function on $[\pi, 2\pi]$.